\documentclass[11 pt]{article}
\usepackage{amsmath,amssymb,amsfonts,amsthm,graphicx,cite,epsfig,makecell, array,mathrsfs}
\usepackage[varg]{pxfonts}
\usepackage{wrapfig}
\usepackage{enumitem}

\usepackage[colorlinks = true,
            linkcolor = red,
            citecolor = blue]{hyperref}
\usepackage{titlesec, float}
\theoremstyle{plain}
\newtheorem{theorem}{Theorem}[section]
\newtheorem{lemma}[theorem]{Lemma}
\theoremstyle{definition}
\newtheorem{definition}[theorem]{Definition}
\newtheorem{corollary}[theorem]{Corollary}

\theoremstyle{remark}
\newtheorem{remark}{\sc \textbf{Remark}}

\newcolumntype{?}{!{\vrule width 1.4pt}}

\def\correspondingauthor{\footnote{Corresponding author.}}
\titlelabel{\thetitle.\quad}
\renewenvironment{abstract}
               {\list{}{\rightmargin\leftmargin}%
                \item[\textbf{\hspace{8.6mm}Abstract ---}]\relax}
               {\endlist}
\DeclareUrlCommand{\url}{%
    \def\UrlLeft##1\UrlRight{\underline{##1}}}

\date{}
\title{Faber polynomial coefficient estimation of subclass of bi-subordinate univalent functions}
\author{S.A. Saleh $^{1}$, Alaa H. El-Qadeem$^{2}$\,,\, Mohamed A. Mamon$^{3,}$\correspondingauthor{}\vspace{0.1in}\\
\footnotesize{$^{1,3}$Department of Mathematics, Faculty of Science, Tanta University, Tanta 31527, Egypt. }  \\
\footnotesize{$^{2}$ Department of Mathematics, Faculty of Science, Zagazig University, Zagazig 44519, Egypt. }  \\
\footnotesize{e-mail: $^1$ Dr.sasaleh@hotmail.com, ~$^2$ ahhassan@science.zu.edu.eg,} \,$^3$ mohamed.saleem@science.tanta.edu.eg}
\begin{document}
\maketitle
\begin{abstract}
In this paper, a comprehensive subclass of bi-univalent functions class are introduced and investigated. Using the Faber polynomials, estimation of the coefficients $|a_n|$ and certain Fekete-Szeg\"{o} inequality of Maclaurin expansion of functions in this subclass are concluded. Finally, some earlier results are pointed out and improved.
\end{abstract}



\section{Introduction}

Let $\mathcal{A}$ denote the class of analytic functions of the form
\begin{equation}\label{eq.1.1.1}
f(z)=z+\sum\limits_{n=2}^{\infty} a_{n}z^{n},
\end{equation}
normalized by the conditions $f(0)=0\ $and $f^{^{\prime}}(0)=1$ defined in the
open unit disk
\begin{equation*}
U=\left\{  z\in \mathbb{C}:\left\vert z\right\vert <1\right\}  .
\end{equation*}

Let $\mathcal{S}$ be the subclass of $\mathcal{A}$ consisting of all
functions of the form (\ref{eq.1.1.1}) which are univalent in $U$. Let $\varphi$ be an analytic
univalent function in $U$ with positive real part and $\varphi(U)$ be symmetric with
respect to the real axis, starlike with respect to $\varphi(0) = 1$ and $\varphi^{'}(0) > 0$. Ma
and Minda \cite{Minda} gave a unified presentation of various subclasses of starlike
and convex functions by introducing the classes $\mathfrak{S^{*}}(\varphi)$ and $\mathfrak{K}(\varphi)$ of functions
$f\in \mathcal{S}$ satisfying $(zf^{'}(z)/f(z)) \prec \varphi(z)$ and $1+(zf^{''}(z)/f^{'}(z))\prec \varphi(z)$ respectively,
which includes several well-known classes as special case. For example, when $\varphi(z)= (1 + Az)/(1 + Bz)$ with a condition $(-1\leq B < A \leq 1)$, the classes $\mathfrak{S^{*}}(\varphi)$ and $\mathfrak{K}(\varphi)$ converted to the class $\mathcal{S^{*}}[A,B]$ and $\mathcal{K}[A,B]$, respectively, introduced by Janowski \cite{Janow}. Although, for a special choose of the value of $A=1-2\beta,~B=-1~(0 \leq \beta < 1)$, the classes $\mathcal{S^{*}}[A,B]$ and $\mathcal{K}[A,B]$ reduced to the classes
$\mathcal{S^{*}}(\beta)$ and $\mathcal{K}(\beta)$, respectively, which are the class of starlike and convex functions of order $\beta$.
For anther choose of the function $\varphi(z)={\left((1+z)/(1-z)\right)}^\alpha$, we obtain the classes $\mathcal{S^{*}_{\alpha}}$ and $\mathcal{K_{\alpha}}$ which are the class of strongly starlike and strongly convex functions of order $\alpha~(0<\alpha\leq1)$.

The Koebe one quarter theorem \cite{Duren} ensures that the image of $U$ under
every univalent function $f\in\mathcal{S}$ contains a disk of radius $\frac
{1}{4}$. Thus every univalent function $f$ has an inverse
$f^{-1}$ satisfying
\begin{equation*}
f^{-1}(f(z))=z,\text{ \ }(z\in U)\text{ and }f(f^{-1}(w))=w\text{
}(\left\vert w\right\vert <r_{0}(f),r_{0}(f)\leq\frac{1}{4}).
\end{equation*}

A function $f\in\mathcal{S}$ is said to be bi-univalent in $U$ if
both $f$ and $f^{-1}$\ are univalent in $U.$ Let $\Sigma$ denote the class of
all bi-univalent functions defined in the unit disk $U$. Since $f\in\Sigma
$\ has the Maclaurin series expansion given by (\ref{eq.1.1.1}), a simple calculation
shows that its inverse $g=f^{-1}$ has the series expansion
\begin{eqnarray*}
g(w)&=&f^{-1}(w)\\
&=& w-a_{2}w^{2}+(2a_{2}^{2}-a_{3})w^{3}-... .
\end{eqnarray*}
Examples of functions in the class $\Sigma$ are
\begin{equation*}
\frac{z}{1-z}, \text{ \ \ } -log\left(1-z\right) \text{ \ \ and \ \ }  \frac{1}{2} log\left(\frac{1+z}{1-z}\right) ,
\end{equation*}
and so on. However, the familiar Koebe function is not a member of $\Sigma$.\\
Other common examples of functions in $\mathcal{S}$ such as
\begin{equation*}
z-\frac{z^2}{2} \text{ \ \ and \ \ } \frac{z}{1-z^2}
\end{equation*}
are also not members of $\Sigma$ (see \cite{Sriv2010}).

Many papers concerning bi-univalent functions have been published recently (for mentioned but a few, \cite{Bulut,Alaa1,Taha,B.Frasin}). A function $f\in \Sigma$ is in the class $\mathcal{S}^{*}_{\Sigma}(\beta)$ of bi-starlike function of order $\beta(0 \leq \beta < 1)$, or $\mathcal{K}_{\Sigma}(\beta)$ of bi-convex function of order $\beta$ if both $f$ and $f^{-1}$ are respectively starlike or convex functions of order $\beta$. For $0< \alpha  \leq 1$,
the function $f\in \Sigma$ is strongly bi-starlike function of order $\alpha$ if both the
functions $f$ and $f^{-1}$ are strongly starlike functions of order $\alpha$. The class
of all such functions is denoted by $\mathcal{S}^{*}_{\Sigma,\alpha}$. These classes were introduced by
Brannan and Taha \cite{Taha}. They obtained estimates on
the initial coefficients $|a_2|$ and $|a_3|$ for functions in these classes. The research into $\Sigma$ was started by Lewin \cite{Lewin}.
He focused on problems connected with coefficients and showed that $|a_2|< 1.51$. Subsequently, Brannan and Clunie \cite{Clunie} conjectured that $|a_2|<\sqrt{2}$. Netanyahu \cite{Netanyahu} concluded that $\max|a_2|=\frac{4}{3}$.

The coefficient estimate problem for each of the following Taylor Maclaurin coefficients$|a_n|$, $n\in \{2,3,\cdots\}$ is presumably still an open problem. This is because the bi-univalency requirement makes the behavior of the coefficients of the function $f$ and $f^{-1}$ unpredictable.
The Faber polynomials play an important role in various areas of mathematical
sciences, especially in geometric function theory. The recent publications \cite{Hamidi,Hamidi3} applying the Faber
polynomial expansions to meromorphic bi-univalent functions motivated us to apply this technique to
classes of analytic bi-univalent functions.
In the literature, there are only a few works determining the general coefficient bounds $|a_n|$ for the
analytic bi-univalent functions given by (\ref{eq.1.1.1}) using Faber polynomial expansions.

In this present work, we use the Faber polynomials in obtaining bounds of Maclaurin coefficients $|a_n|,~n\in \mathbb{N}\ {1}$ and bounds for the Fekete-Szeg\"{o} functional $|a_3 -2a^{2}_{2}|$ of a new defined subclass of $\Sigma$ to generalize some earlier results.
\section{Construction of the subclass $\mathcal{H}_{\Sigma}(\tau,\lambda,\delta;\varphi)$}
\noindent

Throughout this section, let us assume that $\varphi$ be an analytic function with positive real part in the unit disc $U$ satisfying $\varphi(0) = 1,~\varphi^{\prime}(0)> 0 $ and $\varphi(U)$ is symmetric with respect to the real axis. Such a function has a series expansion of the form
\begin{equation}\label{eq.5.1.1}
\varphi(z)=1+B_1z+B_2z^2+B_3z^3+\cdots \text{ \ \ \ \ \ \ }( B_1>0 ).
\end{equation}
where $B_n\in \mathbb{R},$ for all $n=2,3,...$.

Using the Faber polynomial \cite{Airault,Airault2} expansion of the functions $f\in \Sigma$ of the form (\ref{eq.1.1.1}), the inverse function $g=f^{-1}$ may be expressed as
\begin{equation}\label{eq.5.1.2}
g(w)=f^{-1}(w)=w+\sum_{n=2}^{\infty} A_n w^n,
\end{equation}
where
\begin{equation}\label{eq.5.1.3}
A_n=\frac{1}{n} \mathcal{K}_{n-1}^{-n}(a_2,a_3,...,a_n).
\end{equation}
Now, for any $p\in\mathbb{Z}:= \{0,\pm1,\pm2,\cdots\}$, the expansion of $\mathcal{K}_{n}^{p}$ is given by
\begin{equation}\label{eq.5.1.4}
\mathcal{K}_{n}^{p}=p a_{n}+\frac{p!}{(p-2)!2!}D_{n}^{2}+\frac{p!}{(p-3)!3!}D_{n}^{3}+\cdots+\frac{p!}{(p-n)!n!}D_{n}^{n},
\end{equation}
where
\begin{eqnarray}\label{eq.5.1.5}
\nonumber D_{n}^{m}&=& D_{n}^{m}(a_1,a_2,...,a_n),\\
&=&\sum_{n=2}^{\infty}\frac{m!}{\mu_1!\mu_2!\mu_3!\cdots\mu_n!}a_{1}^{\mu_1}a_{2}^{\mu_2}a_{3}^{\mu_3}\cdots a_{n}^{\mu_n},
\end{eqnarray}
while $a_1 =1$ and the sum is taken over all non-negative integers $\mu_1,\mu_2,\mu_3,...,\mu_n$ satisfying
\[
\mu_1+\mu_2+\mu_3+\cdots+\mu_n=m,
\]
\[
\mu_1+2\mu_2+\cdots+n\mu_n=n.
\]
It is observed that
\[
D_{n}^{n}(a_1,a_3,...,a_n)=a_{1}^{n}.
\]
Thus, from equation (\ref{eq.5.1.4}) together with (\ref{eq.5.1.5}) we get an expression of $\mathcal{K}_{n-1}^{-n}$ as
\begin{eqnarray*}
\nonumber \mathcal{K}_{n-1}^{-n}(a_2,a_3,...,a_n)&=& \tfrac{(-n)!}{(-2n+1)!(n-1)!}a_{2}^{n-1}+\tfrac{(-n)!}{(2(-n+1))!(n-3)!}a_{2}^{n-3}a_{3}\\
\nonumber &+&\frac{(-n)!}{(-2n+3))!(n-4)!}a_{2}^{n-4}a_{4}\\
\nonumber &+&\frac{(-n)!}{(2(-n+2))!(n-5)!}a_{2}^{n-5}\left(a_{5}+(-n+2)a_{3}^{2}\right)\\
\nonumber &+&\frac{(-n)!}{(-2n+5))!(n-6)!}a_{2}^{n-6}\left(a_{6}+(-2n+5)a_{3}a_{4}\right)\\
&+&\sum_{j\geq7}a_{2}^{n-j}V_j,
\end{eqnarray*}
where such expressions as $(-n)!$ are to be interpreted by
\[
(-n)!:=\Gamma(1-n)= (-n)(-n-1)(-n-2)\cdots ~~~~~~~~(~n\in \mathbb{N}_{0}:=\mathbb{N}\cup\{0\}~),
\]
and $V_j$ ($7\leq j\leq n$) is a homogeneous polynomial in the variables $a_2, a_3,..., a_n$. In particular, in case of $n=2,3,4$ the expression of $\mathcal{K}_{n-1}^{-n}$ is reduced to
\begin{eqnarray}
\nonumber \mathcal{K}_{1}^{-2}&=& -2a_2,\\
\nonumber \mathcal{K}_{2}^{-3}&=&3(2a_{2}^{2}-a_3),\\
\nonumber \mathcal{K}_{3}^{-4}&=&-4(5a_{2}^{3}-5a_2a_3+a_4).
\end{eqnarray}
\begin{definition}\label{de5.1}
Let $\lambda\geq 1, \tau\in \mathbb{C}^{\ast}=\mathbb{C}-\left\{0\right\}$, $0\leq\delta\leq1$ and $f,g\in \Sigma$ given by (\ref{eq.1.1.1}) and (\ref{eq.5.1.2}) respectively, then $f$ is said to be in the class $\mathcal{H}_{\Sigma}(\tau,\lambda,\delta;\varphi)$ if
\begin{equation}\label{eq.5.1.6}
1+\frac{1}{\tau}\left((1-\lambda)\frac{f(z)}{z}+\lambda f^{'}(z)+\delta zf^{''}(z)-1\right)\prec\varphi(z),
\end{equation}
and
\begin{equation}\label{eq.5.1.7}
1+\frac{1}{\tau}\left((1-\lambda)\frac{g(w)}{w}+\lambda g^{'}(w)+\delta zg^{''}(w)-1\right)\prec\varphi(w),
\end{equation}
where $z,w\in U$ and $\varphi(z)$ is given by (\ref{eq.5.1.1}).
\end{definition}
\begin{remark}
For special choices of the parameters $\lambda, \tau, \delta$ and the function $\varphi(z)$, the class $\mathcal{H}_{\Sigma}(\tau,\lambda,\delta;\varphi)$ reduced to the following subclasses:
\begin{itemize}
  \item [1.] $\mathcal{H}_{\Sigma}\left(\tau,1,\gamma;\varphi\right)=\Sigma(\tau,\gamma,\varphi)$ which introduced by A.E. Tudor \cite{Tudor} and recently studied by H.M. Srivastava and Deepak Bansal \cite{Deep}.
  \item [2.] $\mathcal{H}_{\Sigma}\left(1,1,0;\varphi\right)=\mathcal{H}_\sigma(\varphi)$ which defined and studied by Rosihan M. Ali et al. \cite{Rosihan}.
  \item [3.] $\mathcal{H}_{\Sigma}\left(1,1,\beta;\left(\frac{1+z}{1-z}\right)^{\alpha}\right)=\mathcal{H}_{\Sigma}(\alpha,\beta)$ which introduced by B.A. Frasin \cite{B.Frasin}.
  \item [4.] $\mathcal{H}_{\Sigma}\left(1,1,0;\left(\frac{1+z}{1-z}\right)^{\alpha}\right)=\mathcal{H}_{\Sigma}^{\alpha}$ which introduced by H.M. Srivastava et al. \cite{Sriv2010}.
  \item [5.] $\mathcal{H}_{\Sigma}\left(1,\lambda,0;\left(\frac{1+z}{1-z}\right)^{\alpha}\right)=\mathcal{B}_{\Sigma}(\alpha,\lambda)$ which is introduced by B.A. Frasin and M.K. Aouf \cite{Frasin}, and recently studied by H.M. Srivastava et al. \cite{Eker}.
  \item [6.] $\mathcal{H}_{\Sigma}\left(1-\gamma,1,\beta;\frac{1+z}{1-z}\right)=\mathcal{H}_{\Sigma}(\gamma,\beta)$ which introduced by B.A. Frasin \cite{B.Frasin}.
  \item [7.] $\mathcal{H}_{\Sigma}\left(1-\alpha,\lambda,\delta;\frac{1+z}{1-z}\right)=\mathcal{N}_{\Sigma}(\alpha,\lambda,\delta)$ which introduced by S. Bulut \cite{Bulut}.
  \item [8.] $\mathcal{H}_{\Sigma}\left(1-\beta,1,0;\frac{1+z}{1-z}\right)=\mathcal{H}_{\Sigma}(\beta)$ which introduced by H.M. Srivastava et al. \cite{Sriv2010}.
  \item [9.] $\mathcal{H}_{\Sigma}\left(1-\beta,\lambda,0;\frac{1+z}{1-z}\right)=\mathcal{B}_{\Sigma}(\beta,\lambda)$ which introduced by B.A. Frasin and M.A. Aouf \cite{Frasin} and recently studied by J.M. Jahangiri and S.G. Hamidi \cite{Jah}.
  \item [10.] $\mathcal{H}_{\Sigma}\left(\tau,1,\gamma;\frac{1+Az}{1+Bz}\right)=\mathcal{R}_{\gamma,\sigma}^{\tau}(A,B)$ which introduced by A.E. Tudor \cite{Tudor}.
\end{itemize}
\end{remark}
\begin{lemma} \cite{Motamed} \label{lm5.1}
Let $u(z)$ be analytic function in the unit disc $\mathbb{U}$ with $u(0)=0$ and $|u(z)|<1$ for all $z\in U$ with the power series expansion
\[
u(z)=\sum_{n=1}^{\infty} c_n z^n~,
\]
then $|c_n|\leq 1$ for all $n=1,2,3,...$. Furthermore, $|c_n|=1$ for some $n=1,2,3,...$ if and only if
\[
u(z)=e^{i\theta} z^n \text{, \ \ \ } \theta \in \mathbb{R}.
\]
\end{lemma}
\begin{lemma} \cite{Deniz1} \label{lm5.2}
Let the function $p(z)=1 + \sum_{n=1}^{\infty}p_n z^n$ be so that $\Re(p(z))>0$ for $z\in U$. Then for $-\infty < \alpha < \infty$,
\begin{equation}\label{1812}
 \left|p_2 - \alpha p_1^2\right|\leq \left\{ \begin{array}{cc}
                                                 2-\alpha|p_1|^2 & ;~\alpha< \frac{1}{2} \\
                                                 ~ & ~ \\
                                                 2-(1-\alpha)|p_1|^2 & ;~\alpha\geq \frac{1}{2}
                                               \end{array}\right..
 \end{equation}

Let $\varphi(z)=\sum_{n =1}^{\infty} a_n z^n$ be a Schwarz function so that $|\varphi(z)| < 1$, $z\in U$. Set $p(z)=\frac{1+\varphi(z)}{1-\varphi(z)}$ where $p(z)=1+\sum_{n =1}^{\infty} p_n z^n$ is so that $\Re(p(z))>0$ for $z\in U$. Comparing the corresponding coefficients of powers of $z$ yields $p_1 = 2 \varphi_1 $ and $p_2=2(\varphi_2+\varphi_1^2)$. Now, substituting for $p_1$ and $p_2$ and letting $\eta=1-2\alpha$ in (\ref{1812}), we obtain
\begin{equation}\label{1712}
\left|\varphi_2+\eta \varphi_1^2\right|\leq \left\{ \begin{array}{cc}
                                                 1-(1-\eta)|\varphi_1|^2 & ;~\eta> 0 \\
                                                 ~ & ~ \\
                                                 1-(1+\eta)|\varphi_1|^2 & ;~\eta <0
                                               \end{array}\right..
\end{equation}
\end{lemma}
\subsection{Coefficient bounds of members of $\mathcal{H}_{\Sigma}(\tau,\lambda,\delta;\varphi)$}
\noindent

Unless otherwise mentioned, let us assume in the reminder of this section that $z\in U$, $\lambda\geq 1$, $0\leq \delta\leq 1$ and $\tau\in \mathbb{C}-\left\{0\right\}$.
\begin{theorem}\label{thm5.1}
Let $f$ defined by (\ref{eq.1.1.1}) belong to the class $\mathcal{H}_{\Sigma}(\tau,\lambda,\delta;\varphi)$ and $a_k=0 ~~ (2\leq k\leq n-1)$, then
\begin{equation}\label{eq.5.1.10}
|a_n|\leq \frac{B_{1}|\tau|}{1+(n-1)(\lambda+n\delta)} \text{ \ \ \ \ \ } (~n\geq4~).
\end{equation}
\end{theorem}
\begin{proof} Since $f\in \mathcal{H}_{\Sigma}(\tau,\lambda,\delta;\varphi)$, then we have
\begin{equation}\label{eq.5.1.11}
1+\frac{1}{\tau}\left((1-\lambda)\frac{f(z)}{z}+\lambda f^{'}(z)+\delta zf^{''}(z)-1\right)=1+\sum_{n=2}^{\infty}\left(\frac{1+(n-1)(\lambda+n\delta)}{\tau}\right) a_n z^{n-1},
\end{equation}
and since the inverse map $g=f^{-1}$ represented by (\ref{eq.5.1.2}) also belonging to the same subclass, then
\begin{equation}\label{eq.5.1.12}
1+\frac{1}{\tau}\left((1-\lambda)\frac{g(w)}{w}+\lambda g^{'}(w)+\delta zg^{''}(w)-1\right)=1+\sum_{n=2}^{\infty}\left(\frac{1+(n-1)(\lambda+n\delta)}{\tau}\right) A_n w^{n-1}.
\end{equation}
Now, Since $f,g\in \mathcal{H}_{\Sigma}(\tau,\lambda,\delta;\varphi)$, by the definition \ref{de5.1}, there exist two Schwarz functions $u(z)=\sum_{n=1}^{\infty}c_n z^n$ and $v(w)=\sum_{n=1}^{\infty}d_n w^n$ such that
\begin{equation}\label{eq.5.1.13}
1+\frac{1}{\tau}\left((1-\lambda)\frac{f(z)}{z}+\lambda f^{'}(z)+\delta zf^{''}(z)-1\right)=\varphi(u(z)),
\end{equation}
\begin{equation}\label{eq.5.1.14}
1+\frac{1}{\tau}\left((1-\lambda)\frac{g(w)}{w}+\lambda g^{'}(w)+\delta zg^{''}(w)-1\right)=\varphi(v(w)),
\end{equation}
such that
\begin{equation}\label{eq.5.1.15}
\varphi(u(z))=1-\sum_{n=2}^{\infty}B_1 \mathcal{K}_{n-1}^{-1}(c_1,...,c_{n-1};B_1,...,B_{n-1})z^{n-1},
\end{equation}
\begin{equation}\label{eq.5.1.16}
\varphi(v(w))=1-\sum_{n=2}^{\infty}B_1 \mathcal{K}_{n-1}^{-1}(d_1,...,d_{n-1};B_1,...,B_{n-1})w^{n-1},
\end{equation}
where in general $\mathcal{K}^{p}_{n}=\mathcal{K}^{p}_{n}(\rho_1,...,\rho_n,B_1,...,B_n)$ are defined by
\begin{eqnarray}\label{eq.5.1.17}
\nonumber \mathcal{K}^{p}_{n}&=& \frac{p!}{(p-n)!(n)!}\rho_{1}^{n}\frac{B_n}{B_1}+\frac{p!}{(p-n+1)!(n-2)!}\rho_{1}^{n-2}\rho_2\frac{B_{n-1}}{B_1}\\
\nonumber &+&\frac{p!}{(p-n+2)!(n-3)!}\rho_{1}^{n-3}\rho_{3}\frac{B_{n-2}}{B_1}\\
\nonumber &+&\frac{p!}{(p-n+3)!(n-4)!}\rho_{1}^{n-4}\left(\rho_{4}\frac{B_{n-3}}{B_1}+\frac{p-n+3}{2}\rho_{2}^{2}\frac{B_{n-2}}{B_1}\right)\\
\nonumber &+&\frac{p!}{(p-n+4)!(n-5)!}\rho_{1}^{n-5}\left(\rho_{5}\frac{B_{n-4}}{B_1}+(p-n+4)\rho_{2}\rho_{3}\frac{B_{n-3}}{B_1}\right)\\
&+&\sum_{j\geq6}\rho_{1}^{n-j}X_j,
\end{eqnarray}
where $X_j$ is a homogeneous polynomial of degree $j$ in the variables $\rho_1,\rho_2,...,\rho_n$.\\

Now, comparing the coefficients in both sides of equations (\ref{eq.5.1.13}) and (\ref{eq.5.1.14}) after substituting about $\varphi(u(z))$ and $\varphi(v(w))$ from equations (\ref{eq.5.1.15}) and (\ref{eq.5.1.16}), we have
\begin{equation}\label{eq.5.1.18}
    \frac{1+(n-1)(\lambda+n\delta)}{\tau}a_n=-B_{1}\mathcal{K}_{n-1}^{-1}(c_1,...,c_{n-1};B_1,...,B_{n-1}),
\end{equation}
\begin{equation}\label{eq.5.1.19}
    \frac{1+(n-1)(\lambda+n\delta)}{\tau}A_n=-B_1 \mathcal{K}_{n-1}^{-1}(d_1,...,d_{n-1};B_1,...,B_{n-1}).
\end{equation}
Since $a_k=0 ~(2\leq k\leq n-1)$, then from equation (\ref{eq.5.1.3}) it is easy to conclude
\begin{equation}\label{eq.5.1.20}
A_n=-a_n.
\end{equation}
Therefore, equations (\ref{eq.5.1.18}) and (\ref{eq.5.1.19}) reduced to
\begin{equation}\label{eq.5.1.21}
\frac{1+(n-1)(\lambda+n\delta)}{\tau}a_n=B_{1}c_{n-1},
\end{equation}
\begin{equation}\label{eq.5.1.22}
-\frac{1+(n-1)(\lambda+n\delta)}{\tau}a_n=B_{1}d_{n-1}.
\end{equation}
By subtracting equation (\ref{eq.5.1.22}) from equation (\ref{eq.5.1.21}) obtained
\begin{equation}\label{eq.5.1.23}
a_n=\frac{B_{1}\tau\left(c_{n-1}-d_{n-1}\right)}{2\left(1+(n-1)(\lambda+n\delta)\right)}.
\end{equation}
Applying Lemma \ref{lm5.1} for the coefficients $c_{n-1}$ and $d_{n-1}$ in equation (\ref{eq.5.1.23}) which reduced to the desired estimation. The proof is completed.
\end{proof}

By putting $\tau=1-\alpha (0\leq\alpha<1)$ and $\varphi(z)=\frac{1+z}{1-z}~~(B_1=2)$ in Theorem \ref{thm5.1}, we conclude
\begin{corollary}\label{co5.1}\cite[Theorem 2]{Bulut}
Let $f\in \mathcal{N}_\Sigma(\alpha,\lambda,\delta)$ and $a_k=0\text{ \ }(2\leq k\leq n-1)$, then
\begin{equation*}
|a_n|\leq \frac{2(1-\alpha)}{1+(n-1)(\lambda+n\delta)} \text{ \ \ \ \ \ } (n\geq4).
\end{equation*}
\end{corollary}

Let us put $\lambda=1$ in Corollary \ref{co5.2}, we have
\begin{corollary}\label{co5.2}\cite[Theorem 1]{Eker}
Let us consider $f\in \mathcal{N}_{\Sigma}^{(\alpha,\lambda)}$ and $a_k=0\text{ \ }(2\leq k\leq n-1)$, then
\begin{equation*}
|a_n|\leq \frac{2(1-\alpha)}{n(1+\delta(n-1))} \text{ \ \ \ \ \ } (n\geq4).
\end{equation*}
\end{corollary}

Let us put $\delta=0$ in Corollary \ref{co5.2}, we obtain
\begin{corollary}\label{co5.3}\cite[Theorem 1]{Jah}
If $f\in \mathfrak{D}(\alpha,\lambda)$ and $a_k=0\text{ \ }(2\leq k\leq n-1)$, then
\begin{equation*}
|a_n|\leq \frac{2(1-\alpha)}{1+\lambda(n-1)} \text{ \ \ \ \ \ } (n\geq4).
\end{equation*}
\end{corollary}
\begin{theorem}\label{thm5.2}
Let $f\in\mathcal{H}_{\Sigma}(\tau,\lambda,\delta;\varphi)$ and $ B_1\geq |B_2|$, then
{\large
\begin{equation}\label{eq.5.1.25}
    |a_2| \leq \left\{\begin{array}{cc}
                        \frac{B_1\sqrt{B_1}|\tau|}{\sqrt{B_1^2|\tau|(1+2\lambda+6\delta)+(B_1+B_2)(1+\lambda+2\delta)^2}} & if~~B_2<0, B_1+B_2\leq 0 \\
                        ~&~\\
                        \frac{B_1\sqrt{B_1}|\tau|}{\sqrt{B_1^2|\tau|(1+2\lambda+6\delta)+(B_1-B_2)(1+\lambda+2\delta)^2}} & if~~B_2>0, B_1-B_2\leq 0
                      \end{array}
     \right.,
\end{equation}
\begin{equation}\label{eq.5.1.26}
    |a_3| \leq \left\{\begin{array}{cc}
                        \frac{B_1|\tau|}{1+2\lambda+6\delta} &;~~ B_1>|B_2|\\
                        ~ & ~ \\
                        \frac{|B_2\tau|}{1+2\lambda+6\delta} &;~~ B_1<|B_2|
                      \end{array}
     \right.,
\end{equation}
and
\begin{equation}\label{eq.5.1.27}
    |a_3-2a_{2}^{2}|\leq \left\{\begin{array}{cc}
                        \frac{B_1|\tau|}{1+2\lambda+6\delta} &;~~ B_1>|B_2|\\
                        ~ & ~ \\
                        \frac{|B_2\tau|}{1+2\lambda+6\delta} &;~~ B_1<|B_2|
                      \end{array}
     \right..
\end{equation}}
\end{theorem}
\begin{proof}
Lets us set $n=2, n=3$ in the equations (\ref{eq.5.1.18}) and (\ref{eq.5.1.19}), we deduce
\begin{equation}\label{eq.5.1.28}
\frac{1+\lambda+2\delta}{\tau}a_2=B_{1}c_{1},
\end{equation}
\begin{equation}\label{eq.5.1.29}
-\frac{1+\lambda+2\delta}{\tau}a_2=B_{1}d_{1},
\end{equation}
\begin{equation}\label{eq.5.1.30}
\frac{1+2\lambda+6\delta}{\tau}a_3=B_{1}c_{2}+B_{2}c_{1}^{2},
\end{equation}
and
\begin{equation}\label{eq.5.1.31}
\frac{1+2\lambda+6\delta}{\tau}(2a_{2}^{2}-a_3)=B_{1}d_{2}+B_{2}d_{1}^{2}.
\end{equation}
From equations (\ref{eq.5.1.28}) and (\ref{eq.5.1.29}), we deduce
\begin{equation}\label{eq.1}
    c_1=-d_1,
\end{equation}
and
\begin{equation}\label{eq.2}
    a_2=\frac{B_1c_1\tau}{1+\lambda+2\delta}.
\end{equation}

Now, adding equation (\ref{eq.5.1.30}) to (\ref{eq.5.1.31}) obtains
\begin{equation}\label{eq.5.1.32}
a_{2}^{2}=\tau \left(\frac{(B_{1}(c_{2}+d_{2})+B_{2}(c_{1}^2+d_{1}^2))}{2(1+2\lambda+6\delta)}\right).
\end{equation}
\begin{equation}\label{eq.3}
a_{2}^{2}=\frac{B_1\tau}{2(1+2\lambda+6\delta)} \left[\left(c_{2}+\frac{B_{2}}{B_1}c_{1}^2\right)+\left(d_{2}+\frac{B_{2}}{B_1}d_{1}^2\right)\right].
\end{equation}

Firstly, let $B_2<0(\eta=\frac{B_2}{B_1}<0, ~B_1+B_2\geq 0)$ and applying Lemma \ref{lm5.2} with using equation (\ref{eq.1}), we obtain
\begin{equation}\label{eq.4}
|a_{2}|^{2}\leq\frac{B_1\tau}{1+2\lambda+6\delta} \left[1-\left(\frac{B_1+B_2}{B_1}\right)|c_{1}|^2\right].
\end{equation}

By substituting of $c_1$ from equation (\ref{eq.2}), we conclude
\begin{equation}\label{eq.5}
|a_{2}|^2\leq \frac{|\tau|^2B_{1}^3}{B_1^2|\tau|(1+2\lambda+6\delta)+(B_1+B_2)(1+\lambda+2\delta)^2}.
\end{equation}
Taking the square root of the both side of inequality (\ref{eq.5}), we have
\begin{equation}\label{eq.6}
|a_{2}|\leq \frac{|\tau|B_{1}\sqrt{B_{1}}}{\sqrt{B_1^2|\tau|(1+2\lambda+6\delta)+(B_1+B_2)(1+\lambda+2\delta)^2}}.
\end{equation}

Second, let $B_2>0(\eta=\frac{B_2}{B_1}>0, ~B_1-B_2\geq 0)$ and applying Lemma \ref{lm5.2} with using equation (\ref{eq.1}), then
\begin{equation}\label{eq.7}
a_{2}^{2}\leq\frac{B_1\tau}{1+2\lambda+6\delta} \left[1-\left(\frac{B_1-B_2}{B_1}\right)|c_{1}|^2\right].
\end{equation}
By substituting of $c_1$ from equation (\ref{eq.2}), we conclude
\begin{equation}\label{eq.8}
|a_{2}|^2\leq \frac{|\tau|^2B_{1}^3}{B_1^2|\tau|(1+2\lambda+6\delta)+(B_1-B_2)(1+\lambda+2\delta)^2}.
\end{equation}
Taking the square root of the both side of inequality (\ref{eq.8}), we have
\begin{equation}\label{eq.8}
|a_{2}|\leq \frac{|\tau|B_{1}\sqrt{B_1}}{B_1^2|\tau|(1+2\lambda+6\delta)+(B_1-B_2)(1+\lambda+2\delta)^2}.
\end{equation}
Combining the last inequality with inequality (\ref{eq.6}), we obtain the desired estimate on the coefficient $|a_2|$ which given by (\ref{eq.5.1.25}).

In order to deduce the estimation of $|a_3|$, subtracting equation (\ref{eq.5.1.31}) from (\ref{eq.5.1.30}) with using equation (\ref{eq.2}), obtains
\begin{equation}\label{eq.5.1.35}
a_{3}= a_{2}^{2}+\frac{B_{1}\tau(c_2-d_2)}{2(1+2\lambda+6\delta)}.
\end{equation}
By substituting of $a_2^2$ from equation (\ref{eq.5.1.32}) into (\ref{eq.5.1.35}), we conclude
\begin{equation}\label{eq.5.1.355}
a_{3}=  \frac{\tau (B_{1}c_{2}+B_2c_1^2)}{1+2\lambda+6\delta}.
\end{equation}
Taking the modulus of both sides of equation (\ref{eq.5.1.355}), we get
\begin{equation}\label{eq.5.1.354}
|a_{3}|\leq  \frac{B_{1}|\tau|}{1+2\lambda+6\delta}\left|c_2+\frac{B_2}{B_1}c_1^2 \right|.
\end{equation}

By applying Lemma \ref{lm5.2}, let first $B_2<0(\eta=\frac{B_2}{B_1}<0)$, then
\begin{equation}\label{eq.9}
|a_{3}|\leq\frac{B_{1}|\tau|}{1+2\lambda+6\delta}\left[1-\frac{B_1-B_2}{B_1}|c_1|^2 \right].
\end{equation}
If $B_1-B_2>0$, then we must put $|c_1|$ by its least value $|c_1|=0$. Thus
\begin{equation}\label{eq.10}
|a_{3}|\leq\frac{B_{1}|\tau|}{1+2\lambda+6\delta}.
\end{equation}
If $B_1-B_2<0$, then we must put $|c_1|$ by its maximum value $|c_1|=1$ (using Lemma \ref{lm5.2}). Thus
\begin{equation}\label{eq.11}
|a_{3}|\leq\frac{B_{2}|\tau|}{1+2\lambda+6\delta}.
\end{equation}

Second, let us put $B_2>0(\eta=\frac{B_2}{B_1}>0)$, then
\begin{equation}\label{eq.12}
|a_{3}|\leq\frac{B_{1}|\tau|}{1+2\lambda+6\delta}\left[1-\frac{B_1+B_2}{B_1}|c_1|^2 \right].
\end{equation}
If $B_1+B_2>0$, then we must put $|c_1|$ by its least value $|c_1|=0$. Thus
\begin{equation}\label{eq.13}
|a_{3}|\leq\frac{B_{1}|\tau|}{1+2\lambda+6\delta}.
\end{equation}
If $B_1+B_2<0$, then we must put $|c_1|$ by its maximum value $|c_1|=1$ (using Lemma \ref{lm5.1}). Thus
\begin{equation}\label{eq.14}
|a_{3}|\leq\frac{-B_{2}|\tau|}{1+2\lambda+6\delta}.
\end{equation}
By comparing the estimates of $|a_3|$ in relations from (\ref{eq.10}) to (\ref{eq.13}) which obtain the desired estimate given by (\ref{eq.5.1.26}).
Finally, using equation (\ref{eq.5.1.31}), gives
\begin{equation}\label{eq.5.1.36}
a_3-2a_{2}^{2}=\frac{-\tau( B_{1}d_2+B_2d_1^2)}{1+2\lambda+6\delta}.
\end{equation}
Using the same technique in proving the estimate of $|a_3|$, we get the desired estimate given by (\ref{eq.5.1.27}), then we prefer to omit it.
\end{proof}

In case of $\lambda=1$, Theorem \ref{thm5.2} becomes
\begin{corollary}\label{co5.4}\cite[Theorem 1]{Deep}
Let $f\in \Sigma(\tau,\delta,\varphi)$, then
\begin{equation*}
    |a_2| \leq \left\{\begin{array}{cc}
                        \tfrac{B_1\sqrt{B_1}|\tau|}{\sqrt{3B_1^2|\tau|(1+2\delta)+4(B_1+B_2)(1+\delta)^2}} & B_2<0 \text{ and } B_1+B_2\geq 0\\
                        ~ & ~ \\
                        \tfrac{B_1\sqrt{B_1}|\tau|}{\sqrt{3B_1^2|\tau|(1+2\delta)+4(B_1-B_2)(1+\delta)^2}} & B_2>0 \text{ and } B_1-B_2\geq 0
                      \end{array}
     \right.,
\end{equation*}
\begin{equation*}
    |a_3| \leq \left\{\begin{array}{cc}
                        \frac{B_1|\tau|}{3(1+2\delta)} & B_1>|B_2|\\
                        ~ & ~ \\
                        \frac{|B_2\tau|}{3(1+2\delta)} & B_1<|B_2|
                      \end{array}
     \right..
\end{equation*}
\end{corollary}
Let us put $\varphi(z)=\left(\frac{1+z}{1-z}\right)^\alpha$, $B_1=2\alpha$ and $B_2=2\alpha^2$, and $\tau=1$ in Corollary \ref{co5.4} we have
\begin{corollary}\label{co5.5}\cite[Theorem 2.2]{B.Frasin}
Let $f\in \mathcal{H}_\Sigma(\alpha,\delta)$, then
\begin{equation*}
    |a_2| \leq \frac{2\alpha}{\sqrt{2(2+\alpha)+4\delta(\alpha+\delta-\alpha\delta +2)}},
\end{equation*}
\begin{equation*}
    |a_3| \leq \frac{2\alpha}{3(1+2\delta)}.
\end{equation*}
\end{corollary}
By putting $\tau=1-\gamma$ and $\varphi(z)=\frac{1+z}{1-z},~B_1=B_2=2,$ in Corollary \ref{co5.4}, we obtain
\begin{corollary}\label{co5.6}\cite[Theorem 3.2]{B.Frasin}
Let $f\in \mathcal{H}_\Sigma(\gamma,\delta)$, then
\begin{equation*}
    |a_2| \leq \sqrt{\frac{2(1-\gamma)}{3(1+2\delta)}},
\end{equation*}
\begin{equation*}
    |a_3| \leq \frac{2(1-\gamma)}{3(1+2\delta)}.
\end{equation*}
\end{corollary}

In case of $\tau=1,~\delta=0$ and $\varphi(z)=\left(\frac{1+z}{1-z}\right)^\alpha,~B_1=2\alpha, B_2=2\alpha^2,$ in Theorem \ref{thm5.2}, we have
\begin{corollary}\label{co5.7}\cite[Theorem 2.2]{Frasin}
Let $f\in \mathcal{B}_\Sigma(\alpha,\lambda)$, then
\begin{equation*}
    |a_2| \leq \frac{2\alpha}{\sqrt{(1+\lambda)^2+\alpha(1+2\lambda-\lambda^2)}},
\end{equation*}
\begin{equation*}
    |a_3| \leq \frac{2\alpha}{1+2\lambda}.
\end{equation*}
\end{corollary}

Let us put $\tau=1-\gamma$ and $\varphi(z)=\frac{1+z}{1-z},~B_1=B_2=2,$ in Theorem \ref{thm5.2}, we obtain
\begin{corollary}\label{co5.8}\cite[Theorem 5]{Bulut}
Let $0\leq \alpha<1$ and $f\in \mathcal{N}_\Sigma(\gamma,\lambda,\delta)$, then
\begin{equation*}
|a_2|\leq \sqrt{\frac{2(1-\gamma)}{1+2\lambda+6\delta}},
\end{equation*}
\begin{equation*}
|a_3|\leq \frac{2(1-\gamma)}{1+2\lambda+6\delta},
\end{equation*}
and
\begin{equation*}
    |a_3-2a_{2}^{2}|\leq \frac{2(1-\gamma)}{1+2\lambda+6\delta}.
\end{equation*}
\end{corollary}

By putting $\delta=0$ in Corollary \ref{co5.8}, gets
\begin{corollary}\label{co5.9}\cite[Theorem 3.2]{Frasin}
If $f$ belong to $\mathcal{B}_\Sigma(\gamma,\lambda)$ and $0\leq\gamma<1$, then
\begin{equation*}
|a_2|\leq \sqrt{\frac{2(1-\gamma)}{1+2\lambda}},
\end{equation*}
\begin{equation*}
|a_3|\leq \frac{2(1-\gamma)}{1+2\lambda},
\end{equation*}
\end{corollary}
\begin{remark}
Some results investigated in Corollaries from \ref{co5.4} to \ref{co5.9} represented an improvement of the estimate of $|a_3|$ of the earlier corresponding results.
\end{remark}

\end{document}